\def\update{August 27, 2009} %
  
 \documentclass{article} 
 
\usepackage[applemac]{inputenc}
\usepackage[pdftex]{graphicx}
\usepackage{hyperref,amssymb,amsthm,amsmath}

\def\m#1{$#1$}

\def\M#1{$$#1$$}

\def\bC{\mathbb{C}}

\def\bQ{\mathbb{Q}}
\def\bR{\mathbb{R}}
\def\bZ{\mathbb{Z}}

\def\bx{\mathbf{x}}

\def\uy{\underline{y}}

\numberwithin{equation}{section}

\newtheorem{theorem}[equation]{Theorem}
\newtheorem{proposition}[equation]{Proposition}
\newtheorem{lemma}[equation]{Lemma}

\def\virgule{\raise 2pt \hbox{,}}

\def\cprime{$'$} \def\cprime{$'$} \def\cprime{$'$}

\def\hfl#1#2{\smash{\mathop{\hbox to 10 mm{\rightarrowfill}}
\limits^{\scriptstyle#1}_{\scriptstyle#2}}}

 \def\hlfl#1#2{\smash{\mathop{\hbox to 10 mm{\leftarrowfill}}
\limits^{\scriptstyle#1}_{\scriptstyle#2}}}

\def\grandt{
\hskip -1.15  cm   
{^{^|}}^{\strut} \hskip 1.5  cm   
\hskip -2.95  cm   
{\raise -.3 true cm\hbox{$\vdots$} }\!
 ^{\vrule height .1 pt depth .1pt width 2.6  cm } 
\! {\raise -.3 true cm\hbox{$\vdots$} }
\hskip -2.45  cm} 

\def\flechedte#1#2{\smash{\mathop{\hbox to 8 mm{\rightarrowfill}}
\limits^{\scriptstyle#1}_{\scriptstyle#2}}\!\! \bullet}

\def\flechedroite#1#2{\smash{\mathop{\hbox to 8 mm{\rightarrowfill}}
\limits^{\scriptstyle#1}_{\scriptstyle#2}}}

\def\bC{{\mathbb C}}

\def\bQ{{\mathbb Q}}
\def\Qbar{{\overline{\mathbb Q}}}

\def\bR{{\mathbb R}}

\def\bZ{{\mathbb Z}}

\def\and{\quad\hbox{ and }\quad}

\def\cqfd{\unskip\kern 6pt\penalty 500
\raise -2pt\hbox{\vrule\vbox to 10pt{\hrule width 4pt
\vfill\hrule}\vrule}\par}
\def\virgule{\raise2pt\hbox{\rm ,}}
\def\pointvirgule{\raise2pt\hbox{;}}

\def\adots{\mathinner{\mkern2mu\raise1pt\hbox{.}
\mkern3mu\raise4pt\hbox{.}\mkern1mu\raise7pt\hbox{.}}}

\def\house#1{\setbox1=\hbox{$\,#1\,$}%
\dimen1=\ht1 \advance\dimen1 by 2pt \dimen2=\dp1 \advance\dimen2 by 2pt
\setbox1=\hbox{\vrule height\dimen1 depth\dimen2\box1\vrule}%
\setbox1=\vbox{\hrule\box1}%
\advance\dimen1 by .4pt \ht1=\dimen1
\advance\dimen2 by .4pt \dp1=\dimen2 \box1~ \relax}

\def\grandt{
\hskip -1.15  cm   
{^{^|}}^{\strut} \hskip 1.5  cm   
\hskip -2.95  cm   
{\raise -.3 true cm\hbox{$\vdots$} }\!
 ^{\vrule height .1 pt depth .1pt width 2.6  cm } 
\! {\raise -.3 true cm\hbox{$\vdots$} }
\hskip -2.45  cm}

\def\us{\underline{s}}

\def\uy{\underline{y}}
\def\uz{\underline{z}}

\def\and{\quad\hbox{and}\quad}

\begin{document}

 \def\entete{
 
 \noindent
\null\hskip -1 true cm 
SASTRA
  \hfill
 \\
\null\hskip -1 true cm 
\\
\null\hskip -1 true cm  
 \hfill 
  {\small \it Updated: \rm \update}

  \medskip

 \begin{center}
 
 {\bf

 \bigskip
 
{\large 
  Auxiliary functions in transcendental number theory.}
   
 \smallskip
 \it
 Michel Waldschmidt
 
 }
 
 \bigskip

 \end{center}

   \bigskip
    
    }

     \entete 
    
\tableofcontents

\section*{Acknowledgement}

This survey is based on lectures given for the first time at the International Conference on Number Theory, Mathematical Physics, and Special Functions,    which was organised  at  Kumbakonam (Tamil Nadu, India)  by Krishna Alladi at the  Shanmugha Arts, Science, Technology, Research Academy (SASTRA) in December 2007.   Further lectures on this topic were given by the author at the Arizona Winter School AWS 2008 (Tucson, Arizona, USA), on Special Functions and Transcendence,
organized  in March 2008 by Matt Papanikolas, David Savitt and Dinesh Thakur. Furthermore, this topic was part of a course given by the author at the Fields Institute 
Summer School in Analytic Number Theory and Diophantine Approximation
at the University of Ottawa (Ontario, Canada), organized  in July 2008 by Nathan Ng and Damien Roy. 

The author wishes to thank  Krishna Alladi who suggested that this paper be published in the collection of SASTRA conferences he is editing.  The author is also grateful to Claude Levesque and Paul Voutier for a number of clever remarks on a preliminary version of this survey.

\goodbreak

\section*{Abstract}
We discuss the role of auxiliary functions in the development of transcendental number theory. 
 
 Initially,  auxiliary functions were completely explicit
  (\S~\ref{S:ExplicitFunctions}). 
The earliest transcendence proof is due to Liouville (\S~\ref{SS:AuxiliaryFunctionLiouville}), who produced the first explicit examples of transcendental numbers at a time where their existence was not yet known; in his proof, the auxiliary function is just a polynomial in one variable. Hermite's proof of the transcendence of \m{e} (1873) is much more involved, the auxiliary function he builds 
   (\S~\ref{SS:HermiteAuxiliaryFunction})   is the first example of the Padé approximants
   (\S~\ref{SS:Pade}), which can be viewed as a far reaching generalization of continued fraction expansion   \cite{MR649083,MR1083352}. Hypergeometric functions 
   (\S~\ref{SS:Hypergeometry}) are among the best candidates for using Padé approximations techniques.

Another tool, which plays the role of auxiliary functions, is produced by interpolation formulae   (\S~\ref{S:InterpolationMethods}).
They occurred in the theory after a question by Weierstra\ss\   (\S~\ref{SS:Weierstrass}) on the so-called {\it exceptional set} $S_f$ of a transcendental function $f$, which is the set of algebraic numbers $\alpha$ such that $f(\alpha)$ is algebraic. The answer to his question is that any set of algebraic numbers is the exceptional set of some transcendental function $f$; this  shows that one should add further conditions in order to get transcendence criteria. One way is to replace {\it algebraic number} by {\it rational integer}: this gives rise to the study of integer--valued entire functions 
  (\S~\ref{SS:IVEF}) with the works of  G.~P\'olya (1915),  A.O.~Gel'fond (1929) and many others. 
The connexion with transcendental number theory may not have been clear until the solution by A.O.~Gel'fond in 1929 of the question of the transcendence of $e^\pi$, a special case of Hilbert's seventh problem (\S~\ref{SS:TranscendenceEpuissancePi}). 
 Along these lines, recent developments are due to T.~Rivoal, who renewed forgotten rational interpolation formulae  (1935) of R.~Lagrange
   (\S~\ref{SS:Lagrange interpolation}).

The simple (but powerful) construction by Liouville was extended to several variables by A.~Thue
  (\S~\ref{SSS:OriginThueSiegelLemma}), who introduced the Dirichlet's box principle (pigeonhole principle) 
  (\S~\ref{S:AuxiliaryFnBoxPrinciple})  into the topic of Diophantine approximation in the early 1900's. In the 1920's, Siegel   
  (\S~\ref{SS:ThueSiegelLemma}) developed this idea and applied it in 1932 to transcendental number theory. This gave rise to the Gel'fond--Schneider method   (\S~\ref{SSS:SGS}), which produces the  Schneider--Lang Criterion in one 
  (\S~\ref{SSS:SchneiderLang}) 
  or several 
  (\S~\ref{SSS:SeveralVariables}) variables. Among many developments of this method are results on modular functions
  (\S~\ref{SSS:Modular}).
Variants of the auxiliary functions produced by Dirichlet's Box Principle are universal auxiliary functions, which have small Taylor coefficients at the origin
  (\S~\ref{SS:UniversalAuxiliaryFunction}). 
Another approach, due to K.~Mahler   (\S~\ref{SS:MahlerMethod}), involves auxiliary functions whose existence is deduced from linear algebra instead of Thue--Siegel Lemma \ref{L:TS}.

  In 1991, M.~Laurent introduced interpolation determinants 
    (\S~\ref{S:InterpolationDeterminants}). 
  Two years later, J.B.~Bost used Arakhelov theory      (\S~\ref{S:Arakelov})   to prove slope inequalities, which dispenses of the choice of bases.

  \goodbreak

\section{Explicit functions}\label{S:ExplicitFunctions}

\subsection{Liouville}\label{SS:AuxiliaryFunctionLiouville}

The first examples of transcendental numbers were produced by  Liouville \cite{Liouville1844} in 1844. At that time, it was not yet known that transcendental numbers exist. The idea of Liouville is to show that all algebraic real numbers $\alpha$   are badly approximated by rational numbers. The simplest example is a rational number $\alpha=a/b$: for any rational number $p/q\not=a/b$, the inequality 
$$
\left|\frac{a}{b}- \frac{p}{q}\right|\ge \frac{1}{bq}
$$
holds. For an irrational real number $x$, on the contrary,
for any $\epsilon>0$ there exists a rational number $p/q$ such that
$$
0<\left| x -  \frac{p}{q}\right|\le \frac{\epsilon}{q}\cdotp
$$
This yields an irrationality criterion, which is the basic tool for proving the irrationality of specific numbers: {\it a real number $x$ is irrational if and only if there exists a sequence $(p_n/q_n)_{n\ge 0}$ 
of distinct rational numbers with }
$$
\lim_{n\rightarrow\infty}  q_n \left| x -  \frac{p_n}{q_n}\right| =0.
$$
This criterion is not too demanding: the quality of the approximation is not very strong. Indeed, for a given irrational number $x$,   it is known that there exist much better rational approximations, since 
 there exist infinitely many rational numbers $p/q$ for which
$$
\left| x -  \frac{p}{q}\right|\le \frac{1}{2q^2}\cdotp
$$
It is a remarkable fact that,  on the one hand, there exist such sharp approximations, and on the other hand, we are usually not able to produce (let alone, to show the existence of) much weaker   rational approximations. Also, in examples like $\zeta(3)$, there is essentially a single known explicit sequence of rational approximations which arises from the irrationality proofs -- of course, once the irrationality is established, the existence of other much better approximations follows, but, so far, no one is able to produce an explicit sequence of such approximations.

Liouville extended the irrationality criterion into a transcendence criterion. 
The proof by Liouville involves the irreducible polynomial $f\in\bZ[X]$ of the given irrational algebraic number $\alpha$. Since 
 \m{\alpha} is algebraic, there exists an irreducible polynomial  \m{f\in\bZ[X]} such that \m{f(\alpha)=0}. 
Let \m{d} be the degree of \m{f}.  For $p/q\in\bQ$, the number   \m{q^d f(p/q)} is a non--zero rational integer, hence
\M{
|f(p/q)|\ge \frac{1}{q^d}\cdotp
}
On the other hand, it is easily seen that there exists a constant $c(\alpha)>0$, depending only on $\alpha$ (and its irreducible polynomial $f$),  such that 
\M{
|f(p/q)|\le c(\alpha) \left| \alpha-\frac{p}{q}\right|.
}
An explicit value for a suitable $c(\alpha)$ is given, for instance, in Exercice 3.6 of  \cite{MR1756786}. 
Therefore, 
\M{
 \left| \alpha-\frac{p}{q}\right|\ge \frac{c'(\alpha)}{q^d}\virgule
}
with $c'(\alpha)=1/c(\alpha)$.

Let  \m{\xi} be  a real number such  that, for any  \m{\kappa>0},  there exists a rational number \m{p/q}  with \m{q\ge 2} satisfying 
 \M{
0< \left|\xi-\frac{p}{q}\right|< \frac{1}{q^\kappa}\cdotp
 }
It follows from Liouville's inequality that \m{\xi} is transcendental. Real numbers satisfying this assumption are called {\it Liouville's numbers}. The first  examples  \cite{Liouville1844}, produced by Liouville in 1844,  involved properties of continued fractions, but already in the first part of his 1844   note, he considered series 
$$
\sum_{m\ge 1} \frac{1}{\ell^{m!}}
$$
for $\ell\in\bZ_{\ge 2}$. Seven years later, in \cite{Liouville1851}, he refers to a letter from Goldbach to Euler for numbers 
$$
\sum_{m\ge 1} \frac{k_m}{10^{m!}}
\virgule
$$
where $(k_m)_{m\ge 1}$ is a sequence of integers in the range $\{0,\ldots,9\}$. Next (see p.~140 of  \cite{Liouville1851}), he uses the same argument to prove the irrationality of $$
\sum_{m\ge 1} \frac{1}{\ell^{m^2}} \virgule
$$
whose transcendence has been  proved only in 1996 by Yu.~V.~Nesterenko 
 \cite{MR97m:11102,
MR1837822}. 

We consider below (\S~\ref{SSS:OriginThueSiegelLemma}) extensions of Liouville's result by Thue, Siegel, Roth and Schmidt.

\subsection{Hermite}\label{SS:HermiteAuxiliaryFunction}

During his course at the École Polytechnique in 1815 (see \cite{Stainville}), Fourier gave a simple  proof for the irrationality of $e$, which can be found in many textbooks
(for instance, Th.~47 Chap.~4 \S~7 of Hardy and Wright \cite{MR2445243}). The idea is to  truncate the Taylor expansion at the origin of the exponential function. In this proof, the   auxiliary function is the tail of the Taylor expansion of the exponential function 
$$
e^z-\sum_{n=0}^N \frac{z^n}{n!}  \virgule
$$
which one specializes at $z=1$. As noticed by F.~Beukers, the proof becomes even shorter if one specializes at $z=-1$. 
This proof has been revisited in 1840 by Liouville  \cite{Liouville1840}, who succeeded to extend the argument and to prove that $e^2$ is not a quadratic number. This result is quoted by Hermite in his memoir \cite{JFM05.0248.01}. 
Fourier's argument produces rational approximations to the number $e$, which are sharp enough to prove the irrationality of $e$,  but not the transcendence. The denominators of these approximations are $N!$. One of Hermite's  ideas is to look for other rational approximations with less restrictive restrictions on the denominators. Instead of the auxiliary functions $e^z-A(z)$ for some $A\in\bQ[z]$, Hermite introduces more general auxiliary functions $R(z)=B(z)e^z-A(z)$. He finds a polynomial $B$ such that the Taylor expansion at the origin of $B(z)e^z$ has a large gap: he calls $A(z)$ the polynomial part of the expansion before the gap, so that the auxiliary function $R(z)$ has a zero of high multiplicity at the origin. Hermite gives explicit formulae for $A$,  $B$ and $R$.  In particular, the polynomials $A$ and $B$  have rational coefficients -- the question is homogeneous, one may multiply by a denominator to get integer coefficients. 
Also, he obtains upper bounds for these integer coefficients (they are not too large) and for the modulus of the remainder (which is small on a given disc). As an example, given $r\in\bQ\setminus\{0\}$ and  $\epsilon>0$,  one can use this construction to show the existence of $A$, $B$ and $R$ with $0<|R(r)|<\epsilon$.  Hence $e^r\not\in\bQ$. This gives another proof of Lambert's result on the irrationality of $e^r$ for $r\in\bQ\setminus\{0\}$, and this proof extends to the irrationality of $\pi$ as well \cite{Breusch,NesterenkoIrrationalitePi}.  

Hermite \cite{JFM05.0248.01} 
   goes much  further, since he obtains the transcendence of $e$. 
To achieve this goal,  he considers {\it simultaneous rational approximations to the exponential function}, in analogy with Diophantine approximation. The idea is as follows. 
Let \m{a_0,a_1,\ldots,a_m}  be rational integers and \m{B_0,B_1,\ldots,B_m} be polynomials in \m{\bZ[x]}. For \m{1\le k\le m}, define
\M{
R_k(x)=B_0(x)e^{kx}-B_k(x).
}
Set \m{b_j=B_j(1)}, \m{0\le j\le m} and \M{
R=a_1R_1(1)+\cdots+a_mR_m(1).
}
The numbers $a_j$ and $b_j$ are rational integers, hence 
$$
a_0b_0+a_1b_1+\cdots+a_mb_m=b_0(a_0+a_1e +a_2 e^2+\cdots+a_me^m)-R
$$ 
also. Therefore, if one can prove \m{0<|R|<1}, then one deduces
\M{
a_0+a_1 e+\cdots+a_m e^m\not=0.
}
Hermite's construction is more general: he produces rational  approximations to the functions $1,e^{\alpha_1 x},\ldots,e^{\alpha_m x} $,
when  \m{\alpha_1,\ldots,\alpha_m} are pairwise distinct complex numbers. Let 
\m{n_0,\ldots,n_m}  be rational integers, all \m{\ge 0}. Set  \m{N=n_0+\cdots+n_m}. 
Hermite constructs explicitly polynomials  \m{B_0,\, B_1,\ldots,\, B_m}, with \m{B_j} of degree  \m{N-n_j},  such that each of the functions
\M{
B_0(z) e^{\alpha_k z}-B_k(z),\quad (1\le k\le m)
}
has a zero at the origin of multiplicity at least  \m{N}. 

Such functions are now known as {\it Padé approximations of the second kind} (or {\it of type II}).

\subsection{Padé approximation}\label{SS:Pade}

In his thesis in 1892, 
H.E.~Padé 
 studied systematically the 
approximation of complex
analytic functions by 
rational functions. 
See Brezinski's papers  \cite{MR649083,MR1083352}, where further references to previous works by Jacobi (1845), Fa\'a di Bruno (1859), Sturm, Brioschi, Sylvester, Frobenius (1870), Darboux (1876), Kronecker (1881) are given.

There are two dual points of view, giving rise to the two types of {\it Padé Approximants}
 \cite{MR99a:11088b}. 

Let  \m{f_0,\ldots,f_m} be complex functions which are analytic near the origin  and 
\m{n_0,\ldots,n_m}  be non--negative  rational integers. Set \m{N=n_0+\cdots+n_m}. 

{\it Padé approximants of type II } are 
polynomials  \m{B_0,\, \ldots, B_m} with \m{B_j} having degree  \m{\le N-n_j}, such that each of the functions 
\M{
B_i(z) f_j(z)- B_j(z) f_i(z)\quad (0\le i< j \le m)
}
has a zero at the origin of multiplicity   \m{\ge N+1}.

{\it  Padé approximants of  type I } are  polynomials \m{P_1,\, \ldots, P_m}, with  \m{P_j} of degree  \m{\le n_j}, such that the function 
\M{
P_1(z) f_1(z) +\cdots+P_m(z) f_m (z)
}
has a zero at the origin of multiplicity at least \m{N+m-1}. 

For type I as well as type II, existence of Padé approximants follow from linear algebra: one compares the number of equations which are produced by the vanishing conditions on the one hand, with the number of coefficients of $P$, considered as unknowns, on the other. Unicity of the solution, up to a multiplicative constant (the linear system of equations is homogeneous)   is true only in specific cases ({\it perfect systems}): it amounts to proving that the matrix of the system of equations is regular.

These approximants were also studied by  Ch.~Hermite for the exponentials functions in 1873 and 1893;  later,  in 1917, he gave further integral formulae for the remainder.        For transcendence purposes, Padé approximants of  type I  have been used for the first time  in 1932 by K.~Mahler  \cite{Mahler1931}, who produced effective versions of the transcendence theorems by Hermite, Lindemann and Weierstra\ss.

In the theory of Diophantine approximation, there are transference theorems, initially  due to Khintchine  (see, for instance, 
\cite{MR0349591,LaurentMichel0703146}
). Similar transference properties for Padé approximation have been considered by H.~Jager
\cite{MR0163099} 
 and J.~Coates 
\cite{MR0202665,
MR0218333
}. 

\subsection{Hypergeometric methods}\label{SS:Hypergeometry}

Explicit Padé approximations are known only for restricted classes of functions; however, when they are available, they often  produce very sharp Diophantine estimates. Among the best candidates for having explicit Padé Approximations are the hypergeometric functions. A.~Thue \cite{MR0460050}  developed this idea in the early 20th Century and was able to solve explicitly several classes of Diophantine equations.   
There is a contrast between the measures of irrationality, for instance, which can be obtained by hypergeometric methods, and those produced by other methods, like Baker's method (\S~\ref{SSS:SGS}): typically, hypergeometric methods produce numerical constants with one or two digits (when the expected value is something like $1$ or $2$), where Baker's method produces constants with several hundreds digits. On the other hand,  Baker's method works in much more general situations. Compared with the Thue--Siegel--Roth--Schmidt's method (\S~\ref{SSS:OriginThueSiegelLemma}),  it has the great advantage of being explicit. 

Among many contributors to this topic, we quote 
A.~Thue, C.L.~Siegel, A.~Baker, G.V.~Chudnovskii, M.~Bennett, P.~Voutier, G.~Rhin, C.~Viola, T.~Rivoal\dots
\ 
These works also involve sorts of auxiliary functions (integrals) depending on parameters which need to be suitably selected in order to produce sharp estimates.

 Chapter 2 of 
 \cite{MR99a:11088b} 
 deals with effective constructions in transcendental number theory and includes two sections (\S~6 and \S~7) on generalized hypergeometric functions and series.

\section{Interpolation methods}\label{S:InterpolationMethods}

  We discuss here 
another type of auxiliary function, which occurred in works related with a question of Weierstra\ss\ on the exceptional set of a transcendental entire function. Recall that an {\it entire function} is  a complex valued  function which is analytic in $\bC$. A function $f$ is {\it algebraic} (over $\bC(z)$) if  $f$ is a solution of a functional equation $P(z,f(z))=0$ for some non--zero polynomial $P\in\bC[X,Y]$. An entire function is algebraic  if and only if it is a polynomial. A function which is not algebraic is called {\it transcendental}\footnote{A polynomial whose coefficients are not all algebraic numbers is an algebraic function, namely is  algebraic over $\bC(z)$, but is a transcendental element over $\bQ(z)$. However, as soon as a polynomial assumes algebraic values at infinitely many points (including derivatives),  its coefficients are algebraic. Therefore, for the questions we consider, it makes no difference to consider algebraicity of functions over $\bC$ or over $\bQ$. }.

\subsection{Weierstra\ss\ question}\label{SS:Weierstrass}

Weierstra\ss\ (see \cite{MR58:10772}) initiated the question of investigating the set of algebraic numbers where a given   transcendental 
entire function \m{f} takes algebraic values.

Denote by $\Qbar$ the {\it  field of algebraic numbers } (algebraic closure of $\bQ$ in $\bC$). 
For an entire function $f$, we define the {\it exceptional set $S_f$} of $f$ as the set of algebraic numbers $\alpha$  such that $f(\alpha)$ is algebraic:
$$
S_f:=\bigl\{\alpha\in\Qbar\; ; \; f(\alpha)\in\Qbar\bigr\}.
$$
For instance, Hermite--Lindemann's Theorem on the transcendence of $\log\alpha$ and $e^\beta$ for $\alpha$ and $\beta$ algebraic numbers is the fact that the exceptional set of the function $e^z$ is $\{0\}$. Also, the exceptional set of $e^z+e^{1+z}$ is empty, by the Theorem of Lindemann--Weierstrass. 
The exceptional set of functions like $2^z$ or $e^{i\pi z}$ is $\bQ$, as shown by the Theorem of Gel'fond and Schneider. 

The exceptional set of a polynomial is $\Qbar$ if the polynomial has algebraic coefficients, otherwise it is finite. Also, any finite set of algebraic numbers is the exceptional set of some entire function: for $s\ge 1$ the set $\{\alpha_1,\ldots,\alpha_s\}$ is the exceptional set of the polynomial 
$\pi (z-\alpha_1)\cdots (z-\alpha_s)\in\bC[z]$ and also of the transcendental entire function $(z-\alpha_2)\cdots(z-\alpha_s)e^{z-\alpha_1}$. 
Assuming Schanuel's conjecture,  further explicit examples of exceptional sets for entire functions  can be produced, for instance
$\bZ_{\ge 0}$ or $\bZ$.

The study of exceptional sets started in 1886  with a letter of Weierstrass to Strauss.  This study  was  later developed by   Strauss, Stäckel, Faber -- see 
\cite{MR58:10772}. 
Further results are due to 
van der Poorten, Gramain, Surroca  and others
(see 
\cite{MR89a:11011,SurrocaIMRN}). 

 Among the results which were obtained, a typical one is the following: 
{\it if \m{A} is a countable subset of \m{\bC} and if
\m{E} is a dense subset of \m{\bC}, there exist transcendental entire functions \m{f} mapping \m{A} into \m{E}.}

Also, van der Poorten noticed in \cite{MR372699}  that  there are  transcendental entire functions \m{f} such that \m{D^k f(\alpha)\in\bQ(\alpha)} for all \m{k\ge 0} and all algebraic  \m{\alpha}.

The question of possible sets $S_f$ has been solved in \cite{AWS}: {\it any  set of algebraic numbers is the exceptional set  of some transcendental entire function}.  Also multiplicities can be included, as follows: 
 define the {\it exceptional set with multiplicity} of a transcendental entire  function $f$ as the subset of $(\alpha,t)\in\Qbar\times\bZ_{\ge 0}$ such that $f^{(t)}(\alpha)\in\Qbar$. Here, $f^{(t)}$ stands for the $t$-th derivative of $f$.
 
 Then  {\it any subset of 
 $\Qbar\times\bZ_{\ge 0}$ is the exceptional set with multiplicities of some transcendental entire  function $f$.}
More generally, the main result of  \cite{AWS} is the following:

\begin{quote}
{\it
Let $A$ be a countable subset of $\bC$. For each pair $(\alpha,s)$  with $\alpha\in A$,  and    $s\in\bZ_{\ge 0}$,
  let  $E_{\alpha,s}$ be a dense subset of $\bC$. Then there exists a transcendental entire function $f$ such that 
\begin{equation}\label{Eq:exceptional}
\left(\frac{d}{dz}\right)^s f(\alpha)\in E_{\alpha,s}
\end{equation}
for all $(\alpha,s)\in A\times \bZ_{\ge 0}$.
}
\end{quote} 

One may replace $\bC$ by $\bR$: it means that one may take for the sets $E_{\alpha,s}$ dense subsets of $\bR$,  provided that one requires $A$ to be a countable subset of $\bR$.

The proof is a construction of an interpolation series (see \S~\ref{SS:IVEF})  on a sequence where each $w$ occurs infinitely often. The coefficients 
of the interpolation series are selected recursively  to be sufficiently small (and nonzero), so that the sum $f$ of the series  is a transcendental entire  function. 

This process yields uncountably many such functions. Further, one may also require that they are algebraically independent over $\bC(z)$ together with their derivatives. Furthermore, at the same time, 
one may  request  further restrictions on each of these functions $f$. For instance, given any  transcendental function $g$ with $g(0)\neq 0$, one may require  $|f|_R\le |g|_R$ for all $R\ge 0$.

As a very special case of \ref{Eq:exceptional}  (selecting $A$ to be the set $\Qbar$ of algebraic numbers and each $ E_{\alpha,s}$ to be either $\Qbar$  or its complement in $\bC$), one deduces the existence of uncountably many algebraic independent transcendental entire functions $f$ such that any Taylor coefficient at any algebraic point $\alpha$ takes a prescribed value, either algebraic or transcendental.

\subsection{Integer--valued entire functions}\label{SS:IVEF}

A simple measure for the growth of an entire function $f$   is  the real valued function $R\mapsto |f|_R$, where 
$$
|f|_R=\sup_{|z|=R}|f(z)|.
$$
An entire  function $f$ has an {\it order of growth} $\le \varrho$ if for all $\epsilon>0$ the inequality
$$
 |f|_R\le \exp\bigl( R^{\varrho+\epsilon}\bigr)
$$
holds for sufficiently large $R$.

In 1915, G.~P\'olya  \cite{Polya1915} initiated the study of {\it integer--valued entire functions}; he proved that  {\it  if  $f$ is a transcendental entire function such that  
 \m{f(n)\in\bZ} for all \m{n\in\bZ_{\ge 0}}, then}
\begin{equation}\label{E:Polya}
\limsup_{R\rightarrow\infty}\frac{1}{R}\log  |f|_R>0.
\end{equation}
An example is the function $2^z$ for which the left hand side is  $\log 2$. A stronger version of the fact that 
{\it $2^z$ is the ``smallest'' entire transcendental function mapping the positive integers to rational integers } is the estimate
$$
\limsup_{R\rightarrow\infty}2^{-R} |f|_R\ge 1,
$$
which is valid under the same assumptions as for  (\ref{E:Polya}).

P\'olya's  method involves 
 interpolation series:  given an entire function $f$ and a sequence of complex numbers $(\alpha_n)_{n\ge 1}$, define inductively a sequence  $(f_n)_{n\ge 0}$ of entire functions by $f_0=f$ and, for $n\ge 0$,
 $$
  f_n(z)=f_n(\alpha_{n+1})+(z-\alpha_{n+1})f_{n+1}(z).
 $$
 Define, for $j\ge 0$,
 $$
 P_j(z)= (z-\alpha_1)(z-\alpha_2)\cdots(z-\alpha_j).
 $$
One gets  an expansion
$$
f(z)= A(z)+P_n(z)f_n(z),
$$
where
$$
A=a_0+a_1P_1+ \cdots+
a_{n-1} P_{n-1}\in\bC[z]
\quad\hbox{
and
}\quad 
a_n=f_n(\alpha_{n+1}) \quad(n\ge 0).
$$ 
Conditions for such an expansion to be convergent as $n\rightarrow\infty$ are known --
see, for instance,  \cite{MR28:376}.

Such interpolation series produce formulae for functions with given values at the sequence of points $\alpha_n$; when some of the $\alpha_n$'s are repeated, these formulae involve the successive derivatives of the function at the given point. 
For instance, for  a constant sequence $\alpha_n=z_0$ for all $n\ge 1$, one obtains the Taylor series expansion of $f$ at $z_0$. These formulae have been studied by  I.~Newton and J-L. Lagrange. 


Analytic formulae for the coefficients $a_n$ and  for the remainder $f_n$ follow from Cauchy's residue Theorem. Indeed, let $x$, $z$, $\alpha_1,\ldots,\alpha_n$ be complex numbers with $x\not\in\{z,\; \alpha_1,\ldots,\alpha_n\}$.  Starting from the easy relation
\begin{equation}\label{E:Hermite}
\frac{1}{x-z}=
\frac{1}{x-\alpha_1}+
\frac{z-\alpha_1}{x-\alpha_1}\cdot 
\frac{1}{x-z}\virgule 
\end{equation}
one deduces by induction the next formula due to Hermite:
\begin{align}
\notag
\frac{1}{x-z}
&=
\sum_{j=0}^{n-1} 
\frac{(z-\alpha_1)(z-\alpha_2)\cdots(z-\alpha_j)}
{(x-\alpha_1)(x-\alpha_2)\cdots(x-\alpha_{j+1})}
\\
\notag
&
\qquad\qquad
+
\frac{(z-\alpha_1)(z-\alpha_2)\cdots(z-\alpha_n)}
{(x-\alpha_1)(x-\alpha_2)\cdots(x-\alpha_n)} 
 \cdot 
\frac{1}{x-z}
\cdotp
\end{align}
Let ${\cal D}$ be an open  disc  containing $\alpha_1,\ldots,\alpha_n$, let ${\cal C}$ denote the  circumference of ${\cal D}$, let ${\cal D}'$ be an open disc containing the closure of ${\cal D}$.
Assume
 $f$ is analytic   in ${\cal D}'$. 
Then
$$
a_j= \frac{1}{2i\pi} \int_{\cal C} \frac{f(x)dx }
{(x-\alpha_1)(x-\alpha_2)\cdots(x-\alpha_{j+1})}
\qquad(0\le j\le n-1)
$$
and
$$
f_n(z)= 
\frac{1}{2i\pi} \int_{\cal C}
\frac{f(x) dx}{(x-\alpha_1)(x-\alpha_2)\cdots(x-\alpha_n )(x-z)}\cdotp
$$
P\'olya applies these formulae to prove that if $f$ is an entire function which  does not grow too fast
and satisfies
 \m{f(n)\in\bZ} for  \m{n\in\bZ_{\ge 0}}, then the  coefficients  $a_n$ in the expansion of $f$ at the sequence $(\alpha_n)_{n\ge 1}= \{0,1,2,\ldots\}$ vanish for sufficiently large \m{n}, hence $f$ is a polynomial. 

  Further works on this topic, using a variety of methods, are due to  
G.H.~Hardy, G.~P\'olya, D.~Sato, E.G.~Straus, A.~Selberg, Ch.~Pisot,
F.~Carlson, F.~Gross,\dots \  -- and A.O.~Gel'fond (see \S~\ref{SS:TranscendenceEpuissancePi}).

\subsection{Transcendence of $e^\pi$}\label{SS:TranscendenceEpuissancePi}

P\'olya's study of the growth  of transcendental entire functions taking integral values at the positive rational integers was extended to the Gaussian integers by A.O.~Gel'fond in 1929
\cite{JFM55.0116.01}. By means of 
interpolation series at the points in  \m{\bZ[i]}, he proved that  
{\it if $f$ is 
a transcendental  entire function  which satisfies  \m{f(\alpha)\in\bZ[i]}  for all  \m{\alpha\in\bZ[i]}, then }
\begin{equation}\label{E:Gelfond}
\limsup_{R\rightarrow\infty}\frac{1}{ R^2}\log|f|_R\ge\gamma.
\end{equation}
 
The example of the Weierstra\ss\ sigma function attached to the lattice $\bZ[i]$
$$
\sigma (z)=z\prod_{\omega\in \bZ[i] \setminus\{0\}}
\left(1-\frac{z}{\omega}\right)
e^{\frac{z}{\omega}+\frac{z^2}{2\omega^2}}
$$
 (which is nothing else than the Hadamard canonical product   with  a simple zero at any point of $\bZ[i]$), shows that the constant $\gamma$  in \ref{E:Gelfond} cannot be larger than  $ \pi/2$. Often, for such problems, dealing with a discrete subset of $\bC$, replacing {\it  integer values} by {\it zero values} gives some hint of what should be expected, at least for the order of growth (the exponent $2$ of $R^2$  in the left hand side of formula (\ref{E:Gelfond})), if not for the value of the constant (the number $\gamma$   in the right hand side of formula (\ref{E:Gelfond})). Other examples of Hadamard canonical products are 
 $$
z \prod_{n\ge 1}  \left(
1-\frac{z}{n}\right)e^{z/n}=
-e^{\gamma z}\Gamma(-z)^{-1}
$$ 
for the set    $\bZ_{>0} = \{1,2,\dots\}$ of positive integers and  
$$
z \prod_{n\in\bZ\setminus\{0\}}  \left(
1-\frac{z}{n}\right)e^{z/n}=\pi^{-1}\sin(\pi z)
$$
for the set   $\bZ$ of rational integers.

 The initial admissible 
value computed by A.O.~Gel'fond in 1929 for  $\gamma$ in (\ref{E:Gelfond}) was pretty small, namely \m{\gamma=10^{-45}}. 
It was improved by several mathematicians, including      Fukasawa, Gruman, 
Masser, until 1981, when
F.~Gramain  
  \cite{MR83g:30028}
  reached the value   \m{\gamma=\pi/(2e)}, which is best possible, as shown by  D.W.~Masser 
  \cite{MR81i:10048} one year earlier. See \cite{AblyMzari} for recent results and extensions to number fields.

This work  of Gel'fond's \cite{JFM55.0116.01} turns out to  have fundamental consequences on the development of transcendental number theory, due to its connexion with the number $e^\pi$.  Indeed, the assertion that the number 
\M{e^\pi=23,140\, 692\, 632\, 779\, 269\, 005\, 729\, 086\, 367\, 
\dots}
is irrational is equivalent to saying that the function  \m{e^{\pi z}} cannot take all its values in \m{\bQ(i)} when the argument \m{z}  ranges over  \m{\bZ[i]}. 
By expanding the function  \m{e^{\pi z}} into an interpolation series at the Gaussian integers, Gel'fond was able to prove the transcendence of $e^\pi$. More generally, Gel'fond proved the transcendence of $\alpha^\beta$ for $\alpha$ and $\beta$ algebraic, $\alpha\not=0$, $\alpha\not=1$ and $\beta$  imaginary quadratic. In 1930, Kuzmin extended the proof to the case where $\beta$ is real quadratic, thus proving the transcendence of $2^{\sqrt{2}}$. The same year, Boehle proved that if $\beta $ is algebraic of degree $d\ge 2$, then one at least of the $d-1$ numbers $\alpha^\beta, \alpha^{\beta^2},\ldots, \alpha^{\beta^{d-1}}$ is transcendental (see, for instance, \cite{MR0214551}).  

The next important step came from Siegel's introduction of further ideas in the theory (see \S~\ref{SSS:SGS}).  

We conclude this subsection by noting that our knowledge of the Diophantine properties of the number $e^\pi$ is far from being complete. It was proved recently by Yu.V.~Nesterenko
(see \S~\ref{SSS:Modular})  that the two numbers $\pi$ and $e^\pi$ are algebraically independent, but there is no proof, so far, that $e^\pi$ is not a Liouville number.

\subsection{Lagrange interpolation}\label{SS:Lagrange interpolation}

Newton-Lagrange interpolation  (\S~\ref{SS:IVEF}) of a function yields a series of {\it polynomials}, namely linear combinations of products  $(z-\alpha_1)\cdots(z-\alpha_n)$. Another type of interpolation has been devised  in \cite{JFM61.0321.02}   
by  another Lagrange (René and not Joseph--Louis), in 1935, 
who introduced instead  a series of {\it  rational } fractions. Starting from the formula
\M{
\frac{1}{x-z}=
\frac{\alpha-\beta}{(x-\alpha)(x-\beta)}+
\frac{x-\beta}{x-\alpha}\cdot \frac{z-\alpha}{z-\beta}\cdot 
\frac{1}{x-z}
}
in place of (\ref{E:Hermite}),
 iterating and integrating as in \S~\ref{SS:IVEF}, 
 one deduces  an expansion
\M{
f(z)=\sum_{j=0}^{n-1} b_j 
\frac{(z-\alpha_1)\cdots  (z-\alpha_j)}
{(z-\beta_1)\cdots  (z-\beta_j)}
+R_n(z).
}
This approach has been developed in 2006 \cite{RivoalIJNT}  by  T.~Rivoal,  who applies it to the Hurwitz zeta function
\M{
\zeta(s,z)=\sum_{k=1}^\infty \frac{1}{(k+z)^s}
\quad(s\in\bC,\; \Re e (s) >1,\; z\in\bC).
}
He expands \m{\zeta(2,z)}
as a Lagrange series in 
\M{
\frac{
 z^2(z-1)^2\cdots (z-n+1)^2}{(z+1)^2\cdots (z+n)^2}
\cdotp
} 
He shows that the coefficients of the expansion
belong to  \m{\bQ+\bQ\zeta(3)}. This enables him  to produce a new proof of Apéry's Theorem on the irrationality of  \m{\zeta(3)}. 

Further, he gives a  new proof of the irrationality of  \m{\log 2},  by expanding 
\M{
 \sum_{k=1}^\infty \frac{(-1)^k}{k+z}
}
into a Lagrange interpolation series. Furthermore,  he gives a new proof of the irrationality of  \m{\zeta(2)},  by expanding the function
\M{
 \sum_{k=1}^\infty
 \left( \frac{1}{k}
- 
 \frac{1}{k+z}\right)
}
as a   {Hermite--Lagrange} series in 
\M{
\frac{
\bigl( z(z-1)\cdots (z-n+1) \bigr)^2}
{(z+1)\cdots (z+n)}
\cdotp
} 
It is striking that these constructions yield exactly the same sequences of rational approximations as the one produced by methods which look very much different \cite{MR2111638}.

 Further developments of the interpolation methods should be possible. For instance, 
 Taylor series are the special case of  Hermite's formula with a single point and multiplicities --- they give rise to  {Padé} approximants.    Multiplicities could also be introduced in  
 Lagrange--Rivoal  interpolation.

\section{Auxiliary functions arising from the Dirichlet's box principle}\label{S:AuxiliaryFnBoxPrinciple}

\subsection{Thue--Siegel lemma}\label{SS:ThueSiegelLemma}

Here is a translation of a statement p.~213 of Siegel's paper \cite{56.0180.05}:

\begin{lemma}[Thue--Siegel]\label{L:TS}
Let 
$$
\begin{matrix}
y_1&=&a_{11} x_1+\cdots+a_{1n}x_n\\
&\vdots&\\
y_m&=&a_{m1} x_1+\cdots+a_{mn}x_n\\
\end{matrix}
$$
be $m$ linear forms in $n$ variables with rational integer coefficients. Assume $n>m$. Let $A\in\bZ_{>0}$ be an upper bound for the absolute values of the $mn$ coefficients $a_{kl}$. Then the system of homogeneous linear equations $y_1=0,\; \dots,\; y_m=0$ has a solution in rational integers $x_1,\ldots,x_n$, not all of which are $0$, with absolute values less than $1+(nA)^{m/(n-m)}$.
\end{lemma}
 
The fact that there is a non-trivial solution is a consequence of linear algebra, thanks to the assumption $n>m$. The point here is that Dirichlet's box principe shows the existence of a non-trivial solution satisfying an explicit upper bound. The estimate for that solution is essential in the proofs due to  A.~Thue, C.L.~Siegel, and later A.O.~Gel'fond, Th.~Schneider, A.~Baker, W.M.~Schmidt  and others (however, K.~Mahler devised another method where such estimate is not required, linear algebra suffices -- see \S~\ref{SS:MahlerMethod}). 

The initial proof by Thue and Siegel relied on the box principle. More sophisticated arguments have been introduced by K.~Mahler, using geometry of numbers, and they yield to a number of developments which we do not survey here (see, for instance,  \cite{0421.10019}).

\subsubsection{The origin of the Thue--Siegel Lemma}\label{SSS:OriginThueSiegelLemma}

The first improvement of Liouville's inequality was reached by A.~Thue in 1909  \cite{JFM40.0265.01,MR0460050}.   Instead of   evaluating the values at \m{p/q} of a polynomial in a single variable (viz{.} the irreducible polynomial of the given algebraic number \m{\alpha}), he considers two approximations \m{p_1/q_1} and   \m{p_2/q_2} of \m{\alpha} and evaluates at the point \m{(p_1/q_1,p_2/q_2)}  a polynomial \m{P} in two variables. This polynomial \m{P\in\bZ[X,Y]} is constructed (or, rather, is shown to exist)
 by means of {Dirichlet}'s box principle (Lemma \ref{L:TS}). 
 The required conditions are that \m{P} has zeroes of sufficiently large multiplicity at \m{(0,0)} and at  \m{(p_1/q_1,p_2/q_2)}. The multiplicity is weighted: this is  what Thue called the {\it index} of \m{P} at a point.   The estimate for the coefficients of the solution of the system of linear equations in Lemma \ref{L:TS}\ plays an important role in the proof. 
 
One of the main difficulties that Thue had to overcome was to produce a {\it zero estimate}, in order to find a non--zero value of some derivative of $P$.

A crucial feature of Thue's argument is that he needs to select a second approximation \m{p_2/q_2} depending on a first one  \m{p_1/q_1}.  Hence, a first very good approximation  \m{p_1/q_1} is required to produce an effective result from this method.  In general, such arguments lead to  sharp estimates for all  \m{p/q} with at most one exception. This approach has been worked out  by J.W.S.~Cassels, H.~Davenport and others to deduce  {\it upper bounds for the number of solutions of certain Diophantine equations}.  However, these results are not {\it effective}, meaning that they do not yield complete solutions of these equations. 
More recently, E.~Bombieri has produced examples where a sufficiently good approximation exists for the method to work in an effective way. Later, he produced {\it effective refinements } to Liouville's inequality by extending the argument (see  
 \cite{MR99a:11088b} Chap.~1 \S~5.4). 

Further improvement of Thue's method were obtained by {C.L.~Siegel } in the 1920's: he developed Thue's method and succeeded in refining his irrationality measure for algebraic real numbers. In 1929,  Siegel  \cite{56.0180.05}, thanks to a further  sharpening  of his previous estimate,  derived his well known theorem on integer points on curves: \emph{the set of integral points on a curve of genus $\ge 1$ is finite}. 

The introduction of the  fundamental memoir 
\cite{56.0180.05} 
 of C.L.~Siegel  in 1929 
stresses the importance of Thue's idea involving the pigeonhole principle. In the second part of this  paper, he extends the Lindemann--Weierstra\ss\ Theorem (on the algebraic independence of $e^{\beta_1},\ldots,e^{\beta_n}$ when $\beta_1,\ldots,\beta_n$ are $\bQ$-linearly independent algebraic numbers) from the usual exponential function to a wide class of entire functions, which he calls {\it $E$-functions}. He also introduces the class of $G$-functions, which has been extensively studied since 1929. 
See also his monograph in  1949 \cite{0039.04402}, 
Shidlovskii's book \cite{MR1033015} 
and the Encyclopaedia volume by  Feldman and Nesterenko
 \cite{MR99a:11088b} %
for $E$--functions,
André's book \cite{MR990016}
and
 \cite{MR99a:11088b} Chap.~5 \S~7 
  for $G$--functions. 
Among many developments related to $G$ functions are the works of Th.~Schneider,  V.G.~Sprindzuck and P.~Dèbes related to algebraic functions (see
 \cite{MR99a:11088b} Chap.~5 \S~7). 

The work of Thue and Siegel  on the  rational approximation to algebraic numbers was extended by many a mathematician,  including Th.~Schneider, A.O.~Gel'fond, 
F.~Dyson,  until 
K.~F.~Roth obtained in 1955 a result which is essentially optimal. In his proof, he  introduces polynomials in many variables.
   
An incredibly  powerful higher dimensional generalization of Thue--Siegel--Roth's Theorem, the  {\it Subspace Theorem}, was obtained in 1970 by  W.M.~Schmidt 
\cite{0421.10019}; 
see also \cite{MR2487729,%
1115.11034}.  
Again, the proof involves a construction of an auxiliary polynomial in several variables, and one of the most difficult auxiliary results is a zero estimate ({\it index theorem}).     

Schmidt's Subspace Theorem,  together with its variants (including effective estimates for the exceptional subspaces, as well as results involving several valuations), have a large number of applications: to Diophantine approximation and Diophantine equations, to transcendence and algebraic independence, to the complexity of algebraic numbers -- see 
\cite{MR2487729}. 

We give here only a simplified statement of this fundamental result, which is already quite deep and very powerful.

\begin{theorem}[Schmidt's Subspace Theorem -- simplified form]\label{Th:SchmidtSubspaceThm}
 For $ m\ge 2$ let $L_1,\ldots,L_{m}$ be independent  linear forms in $m$ variables with algebraic coefficients. Let $\epsilon>0$. Then the set
   $$
   \{\bx =(x_1,\ldots,x_m)\in\bZ^m\; ;\;
   |L_1(\bx)\cdots L_{m}(\bx)|\le |\bx|^{-\epsilon}
   \}
   $$
   is contained in the union of finitely many proper subspaces of $\bQ^m$.

\end{theorem}

\subsubsection{Siegel, Gel'fond, Schneider}\label{SSS:SGS}

In 1932, 
C.L.~Siegel 
 \cite{58.0395.01} 
 obtained the first results on the transcendence of elliptic integrals of the first kind ({\it a Weierstrass  elliptic function cannot have simultaneously algebraic periods and algebraic invariants $g_2$, $g_3$}), by means of a very ingenious argument, which involved an auxiliary function whose existence follows from the Dirichlet's box principle.  This idea turned out to be crucial in the development of transcendental number  theory.

The seventh of the 23 problems raised by 
D.~Hilbert in 1900 is to prove  the  transcendence of the numbers $\alpha^\beta$ for $\alpha$ and $\beta$ algebraic ($\alpha\not=0$, $\alpha\not=1$, $\beta\not\in\bQ$). 
In this statement, $\alpha^\beta$ stands for $\exp(\beta\log\alpha)$, where $\log\alpha$ is any\footnote{The assumption $\alpha\not=1$ can be replaced by the weaker assumption $\log\alpha\not=0$. That means that one can take $\alpha=1$, provided that we select for $\log\alpha$ a non--zero multiple of $2i\pi$. The result allowing $\alpha=1$ is not more general: it amounts to the same to take $\alpha=-1$, provided that one replaces $\beta$ by $2\beta$.} 
 logarithm of $\alpha$.
 The solution was achieved independently by 
A.O.~Gel'fond 
 \cite{JFM60.0163.04} 
 and Th.~Schneider \cite{0010.10501}  in 1934. Consequences, already quoted by Hilbert, are the facts that   $2^{\sqrt{2}}$ and $e^\pi$ are transcendental.  The question of the arithmetic nature of $2^{\sqrt{2}}$ was considered by L.~Euler in 1748 \cite{MR0016336}.

The proofs by Gel'fond and Schneider are different, but both of them rest on some auxiliary function, which arises from Dirichlet's box principle, following Siegel's contribution to the theory. 

Here are the basic ideas behind the methods of Gel'fond and Schneider. 
Let us argue by contradiction and assume that $\alpha$, $\beta$ and $\alpha^\beta$ are all algebraic, with 
$\alpha\not=0$, $\alpha\not=1$, $\beta\not\in\bQ$. 
 Define $K=\bQ(\alpha,\beta,\alpha^\beta)$. By assumption, $K$ is a number field.

A.O.~Gel'fond's proof 
 \cite{JFM60.0163.04} 
 rests on the fact that 
the two entire  functions \m{e^z} and \m{e^{\beta z}} are algebraically independent, they satisfy differential equations with algebraic coefficients and they take simultaneously values in \m{K} for infinitely many \m{z}, viz{.}  \m{ z\in\bZ\log \alpha}.

Th.~Schneider's proof  \cite{0010.10501}   is different: he notices that the two entire functions \m{z} and \m{\alpha^z=e^{z\log\alpha}} are algebraically independent,   they take simultaneously values in \m{K} for infinitely many \m{z}, viz{.} \m{ z\in\bZ+\bZ\beta}. 
He makes no  use of differential equations, since the coefficient $\log\alpha$,  which occurs by derivating the  function $\alpha^z$,  is not algebraic.

Schneider  introduces a polynomial $A(X,Y)\in\bZ[X,Y]$  in two variables and considers the auxiliary function 
$$
F(z)=A(z,\alpha^z)
$$
at the points $m+n\beta$: these values $\gamma_{mn}$ are in the number field $K$.

Gel'fond also 
introduces   a polynomial $A(X,Y)\in\bZ[X,Y]$  in two variables and considers the auxiliary function 
$$
F(z)=A(e^z, e^{\beta z});
$$
the values $\gamma_{mn}$ at the points $m\log\alpha$ of the derivatives  $F^{(n)}(z)$  are again in the number field $K$.

With these notations,  the proofs are similar: the first step is the existence of a non--zero polynomial $A$, of partial degrees bounded by $L_1$ and $L_2$, say,  such that the associated numbers $\gamma_{mn}$ vanish for certain values of $m$ and $n$, say $0\le m<M$, $0\le n<N$. This amounts to showing that a system of linear homogeneous equations has a non--trivial solution;  linear algebra suffices for the existence. In this system of equations, the coefficients are algebraic numbers in the number field $K$, the unknowns are the coefficients of the polynomial $A$.
There are several options at this stage: one may either require  only that the coefficients of $A$ lie  in the ring of integers of $K$, in which case the assumption $L_1L_2>MN $ suffices. An alternative way  is to require the coefficients of $A$ to be in $\bZ$, in which case one needs to assume $L_1L_2>MN[K:\bQ] $. 

This approach is not quite  sufficient for the next steps: one will need estimates for the coefficients of this auxiliary polynomial $A$.  This is where the Thue--Siegel Lemma \ref{L:TS} comes into the picture: by assuming that the number of unknowns, namely $L_1L_2$, is slightly larger than the number of equations, say twice as large, this lemma produces a bound,  for a non--trivial solution of the homogeneous linear system,
which is sharp enough for the rest of the proof. 

The second step is an induction: one proves that $\gamma_{mn}$ vanishes for further values of $(m,n)$. Since there are two parameters $(m,n)$, there are several options for this extrapolation (increasing $m$, or $n$, or $m+n$, for instance), but, anyway, the idea is that if $F$ has sufficiently many zeroes, then $F$ takes rather small values on some disc (Schwarz Lemma), and so do  its derivatives (Cauchy's inequalities). Further, an element of $K$ which is sufficiently small should vanish (by a Liouville type inequality,  or a so-called {\it size } inequality, or else the product formula -- see, for instance,  \cite{MR35:5397,1115.11034,MR99a:11088b,MR1756786} among many references on this topic).

For the last step, there are also several options: one may perform the induction with infinitely many steps and use an asymptotic zero estimate, or else stop after a finite number of steps and prove that some determinant does not vanish. The second method is more difficult and this is the one Schneider succeeded to complete, but his proof can be simplified by pursuing the induction forever.

There is a duality between the two methods.  In  Gel'fond 's proof, replace $L_1$ and $L_2$ by $S_1$ and $S_2$, and replace $M$ and $N$ by $T_0$ and $T_1$; hence,
the numbers $\gamma_{mn}$ which arise are
\M{
\left(\frac{d}{dz}\right)^{t_0} \bigl(e^{(s_1+s_2\beta) z}\bigr) _{z=t_1\log\alpha}.
}
In Schneider's proof,
 replace $L_1$ and $L_2$ by $T_1$ and $T_2$, and replace $M$ and $N$ by $S_0$ and $S_1$; then 
the numbers  $\gamma_{mn}$  which arise are
\M{
\bigl(z^{t_0}\alpha^{t_1z}\bigr)_{z=s_1+s_2\beta}.
}
It is easily seen that the numbers $\gamma_{mn}$ arising from Gel'fond's and Schneider's methods are the same, namely
\begin{equation}\label{Eq:DualiteGS}
 (s_1+s_2\beta)^{t_0}\alpha^{t_1s_1}(\alpha^\beta)^{t_1s_2}.
\end{equation}
See \cite{MR1243105}   
and \S~13.7 of \cite{MR1756786}. 

Gel'fond--Schneider Theorem was extended  in 1966  by A.~Baker \cite{pre05232987}, who proved the more general result that {\it if $\log\alpha_1,\ldots,\log\alpha_n$ are $\bQ$--linearly independent logarithms of algebraic numbers, then the numbers  $1,\log\alpha_1,\ldots,\log\alpha_n$ are linearly independent over $\Qbar$}. 
The auxiliary function used by Baker may be considered as a function of several variables, or as a function of a single complex variable, depending on the point of view
(cf.~\cite{MR50:12931}). 
The analytic estimate (Schwarz lemma) involves merely a single variable. 
The differential equations can be written with a single variable with transcendental coefficients. By introducing several variables, only algebraic coefficients occur. 
See also \S~\ref{SSS:SeveralVariables}.

To be more precise, assume that  $\alpha_1,\ldots,\alpha_n,\alpha_{n+1}$, $\beta_0,\ldots,\beta_n$ are algebraic numbers  which  satisfy
$$
\beta_0+ \beta_1\log\alpha_1+ \beta_n\log\alpha_n=\log\alpha_{n+1}
$$
for some specified values of the logarithms of the $\alpha_j$. Then the  $n+2$ functions of $n+1$ variables
$$
z_0, e^{z_1},\ldots, e^{z_n}, e^{\beta_0z_0+\beta_1z_1+\cdots+\beta_nz_n}
$$
satisfy differential equations with algebraic coefficients and  take algebraic values at the
 integral multiples of the point
$$
(1,\log\alpha_1,\ldots,\log\alpha_n)\in\bC^{n+1}.
$$
This situation is therefore an extension of the setup in Gel'fond's solution of Hilbert's seventh problem, and Baker's method can be viewed as an extension of Gel'fond's method. The fact that all points are on a complex line $\bC(1,\log\alpha_1,\ldots,\log\alpha_n)\subset \bC^{n+1}$ means that Baker's method requires only tools from the theory of one complex variable. 

On the other hand, the corresponding extension of Schneider's method requires several variables: under the same assumptions, consider the functions 
$$
z_0, z_1,\ldots,z_n,e^{z_0}\alpha_1^{z_1}\cdots\alpha_n^{z_n}
$$
and the points in the subgroup of $\bC^{n+1}$ generated by
$$
\bigl(
\{0\}\times \bZ^n
\bigr)
+\bZ(\beta_0,\beta_1,\ldots,\beta_n).
$$
Since Baker's Theorem includes the transcendence of $e$, there is no hope to prove it without introducing the differential equation of the exponential function --  in the factor $e^{z_0}$ of the last function $e^{z_0}\alpha_1^{z_1}\cdots\alpha_n^{z_n}$, the number  $e$ cannot be replaced by an algebraic number! -- hence we need to take also derivatives with respect to $z_0$. For this method, we refer to 
 \cite{MR1756786}. 
 The duality between Baker's method and Schneider's method in several variables is explained below in 
\S~\ref{SS:UniversalAuxiliaryFunction}.

\subsubsection{Schneider--Lang Criterion}\label{SSS:SchneiderLang}

In 1949,
 \cite{MR11:160a} 
 Th.~Schneider produced a very general statement on algebraic values of analytic functions, which can be used as {\it a principle for proofs of transcendence}. This statement includes a large number of previously known results, like the Hermite--Lindemann and Gel'fond--Schneider Theorems. It 
  also contains the so--called {\it Six exponentials Theorem}   \ref{T:SixExp}  
  (which was not explicitly in the literature then). To a certain extent, such statements provide  partial answers to Weierstra\ss\ question (see  
\S~\ref{SS:Weierstrass}) that exceptional sets of transcendental functions are not too large; here, one puts  restrictions on the functions, while in P\'olya's work concerning integer--valued entire functions, the assumptions were mainly on the points and  the values (the mere condition on the functions were that they have a finite order of growth). 

A few years later, in his book 
 \cite{MR19:252f} 
 on transcendental numbers, 
 Schneider gave variants of this statement, which lose some generality but gained in simplicity. 
 
Further simplifications were introduced by S.~Lang in 1964
and the statement which is reproduced in his book on transcendental numbers 
\cite{MR35:5397} 
is the so-called {\it Criterion of Schneider--Lang}
(see also the appendix of  
  \cite{MR2003e:00003} 
as well as
 \cite{MR50:12931} Th.~3.3.1). 

\begin{theorem}\label{T:CriterionSL}

 Let $K$ be a number field and $f_1,\ldots,f_d$ be entire functions in $\bC$. Assume that $f_1$ and $f_2$ are algebraically independent over $K$ and have  finite order of growth. Assume also that they satisfy differential equations: for $1\le i\le d$, assume that the derivative $f'_i$ of $f_i$ is a polynomial in $f_1,\ldots,f_d$ with coefficients in $K$. Then the set $S$ of $w\in\bC$ such that all $f_i(w) $ are in $K$ is finite\footnote{For simplicity, we consider only {\it entire} functions; Schneider--Lang's Theorem deals, more generally, with  meromorphic functions, and this is important for applications,  for instance to elliptic functions. To deal with functions which are analytic in a disc only is also an interesting issue
 \cite{MR784068,MR802729,MR873880,MR948473,MR1145499}.}. 

\end{theorem}

This   statement includes the Hermite--Lindemann Theorem on the transcendence of $e^\alpha$: take 
$$
K=\bQ(\alpha,e^\alpha), \quad
f_1(z)=z, \quad
f_2(z)=e^z, \quad
S=\{ m\alpha\; ;\; m\in\bZ\},
$$ 
as well as the Gel'fond --Schneider Theorem  on the transcendence of $\alpha^\beta$ following Gel'fond's method: take 
$$
K=\bQ(\alpha,\beta, \alpha^\beta),
\quad
f_1(z)=e^z,
\quad
f_2(z)=e^{\beta z},
\quad
S=\{ m\log \alpha\; ;\; m\in\bZ\}.
$$
This criterion \ref{T:CriterionSL}  does not include some of the results which are proved by means of Schneider's method (for instance, it does not contain the {\it Six Exponentials Theorem} \ref{T:SixExp}),
 but there are different criteria (not involving differential equations) for that purpose (see, for instance, \cite{MR35:5397,
MR50:12931,
MR1756786}). 

   Here is the idea of the proof of the Schneider--Lang Criterion  \ref{T:CriterionSL}. 
We argue by contradiction: assume \m{f_1} and \m{f_2} take simultaneously their values in the number field  \m{K} for different values  \m{w_1,\ldots,w_S\in \bC}, where $S$ is sufficiently large. We want to show that there exists a non--zero polynomial \m{P\in \bZ [X_1,X_2]}, such that the function \m{P(f_1,f_2)} is the zero function: this will contradict the assumption that $f_1$ and $f_2$  are algebraically independent.

The first step is to show that there exists a non--zero polynomial \m{P\in \bZ[X_1,X_2]}, such that the function \m{F=P(f_1,f_2)} has a zero of high multiplicity, say $\ge T$, at each \m{w_s}, $(1\le s\le S)$: we consider the system of $ST$ homogeneous linear equations 
\begin{equation}\label{E:SchneiderLang}
\left(\frac{d}{dz}\right)^{t}F(w_s)=0\quad
\hbox{ for }\quad 1\le s\le S
\quad
\hbox{ and }\quad
0\le t<T,
\end{equation}
where the unknowns are the coefficients of $P$.  If we require that the partial degrees of $P$ are strictly less than $L_1$ and $L_2$ respectively, then  the number of unknowns is $L_1L_2$. Since we are looking for a polynomial $P$ with rational integer  coefficients, we need to introduce the degree $[K:\bQ]$ of the number field $K$. 
 As soon $L_1L_2>TS[K:\bQ]$,   there is a non--trivial solution.  Further, the Thue--Siegel Lemma   \ref{L:TS}  produces an upper bound for  the coefficients of \m{P}. This bound is sharp enough if one assumes for instance  $L_1L_2\ge 2TS[K:\bQ]$.
 
 The next step is an induction: the  goal is to prove that \m{F=0}.
By induction on \m{T'\ge T}, one first  proves 
\M{
\left(\frac{d}{dz}\right)^{t}F(w_s)=0\quad
\hbox{ for }\quad 1\le s\le S
\quad
\hbox{ and }\quad
0\le t<T'.
} 
One already knows that these conditions hold for $T'=T$ by (\ref{E:SchneiderLang}).

The proof of this induction is the same as the ones by Gel'fond and Schneider of the transcendence of $\alpha^\beta$,  (see \S~\ref{SSS:SGS}), combining an analytic upper bound (Schwarz Lemma) and an arithmetic lower bound (Liouville's inequality). 
At the end of the induction, one deduces \m{F=0}, which is the desired contradiction with the algebraic independence of \m{f_1} and \m{f_2}.

As we have seen in \S~\ref{SSS:SGS}, such a scheme of proof is characteristic of the Gel'fond--Schneider's method. 

The main analytic argument is Schwarz Lemma for functions of one variable, which produces an upper bound for the modulus of an analytic function having many zeroes. One also uses Cauchy's inequalities in order to bound the moduli of the derivatives of the auxiliary function. 

In this context, a well known open problem raised by Th. Schneider (this is the second in the list of his 8 problems from his book \cite{MR19:252f})   is related with his proof of  the transcendence of $j(\tau)$. Here, $j$ denotes the modular function defined in the upper half plane $\Im{m}(z)>0$, while  $\tau$ is an algebraic point in this upper half plane which is not imaginary quadratic. Schneider himself proved the transcendence of $j(\tau)$, but his proof is not direct, it rests on the use of elliptic functions (one may apply the Schneider--Lang Criterion for meromorphic functions). His  question is to prove the same result by using modular functions. In spite of recent progress on transcendence of values of modular functions (see \S~\ref{SSS:Modular}), this problem is still open. The difficulty lies in the analytic estimate and the absence of a suitable Schwarz  Lemma -- the best  results on this topic  are due to I.Wakabayashi
\cite{MR784068,MR802729,MR873880,MR948473,MR1145499}.

\subsubsection{Higher dimension:   several variables}\label{SSS:SeveralVariables}

In 1941, Th.~Schneider 
 \cite{MR3:266b} 
obtained an outstanding result on the values of Euler's Gamma and Beta functions:
{\it  for any rational numbers $a$ and $b$ such that none of $a$, $b$ and $a+b$ is an integer, the number 
$$
B(a,b)=\frac{\Gamma(a)\Gamma(b)}{\Gamma(a+b)}
$$ 
is transcendental. }

His proof involves a generalization of Gel'fond's method to several variables and yields  a general transcendence criterion for functions satisfying differential equations with algebraic coefficients and taking algebraic values at the points of a large Cartesian product. He applies this criterion to the theta functions associated with the Jacobian of the Fermat curves. His transcendence results apply, more generally, to yield transcendence results on periods of abelian varieties. 

After a suggestion of P.~Cartier, S.~Lang extended the classical results on the transcendence of the values of the classical exponential function to the exponential function of commutative algebraic groups. During this process, he generalized the one dimensional Schneider--Lang criterion  \ref{T:CriterionSL}   to several variables
\cite{MR35:5397} Chap.~IV. 
 In  this higher dimensional criterion, the conclusion is that the set of exceptional values in $\bC^n$ cannot contain a large Cartesian product. This was the generalization to several variables of the fact that in the one dimensional case the exceptional set is finite (with a bound for the number of elements).

A simplified version of the Schneider--Lang Criterion in several variables for Cartesian products is the following 
(\cite{MR1756786} Theorem 4.1).

\begin{theorem}[Schneider--Lang Criterion in several variables] \label{T:SLseveralVariables}
Let
$d$ and $n$ be two integers with $d>n\ge 1$,
$K$ be a number field, and $f_1,\ldots,f_d$ be
algebraically  independent entire functions of finite order
of growth. Assume, for $1\le \nu\le n$ and $1\le i\le d$,
that the partial derivative
$(\partial/\partial z_\nu)f_i$ of $f_i$ belongs to the ring
$K[f_1,\ldots,f_d]$. Further, let $(\uy_1,\ldots,\uy_n)$ be 
a basis of $\bC^n$ over $\bC$. Then the numbers
$$
f_i(s_1\uy_1+\cdots+s_n\uy_n),\qquad\bigl(1\le i\le d,\;
(s_1,\ldots,s_n)\in\bZ^n\bigr) 
$$
do not all belong to $K$.

\end{theorem}

Besides the corollaries already derived by Schneider in 1941, Lang gave further consequences of this result to commutative algebraic groups, especially abelian varieties \cite{MR35:5397} Chap.~IV 
and 
\cite{MR88i:11047}  Chap.~5. 

It is interesting, from an historical point of view, to notice that Bertrand and Masser 
\cite{MR81e:10032} 
succeeded in 1980  to deduce Baker's Theorem form the Schneider--Lang Criterion \ref{T:SLseveralVariables}
(see  
 \cite{MR1756786} \S~4.2). 
 They could also prove the elliptic analog of Baker's result and obtain the linear independence, over the field of algebraic numbers, of elliptic logarithms of algebraic points -- at that time such a result was available only in the case of complex multiplication (by previous work of D.W.~Masser). 

According to 
\cite{MR35:5397}, Historical Note of Chap.~IV, 
M.~Nagata suggested that, in the higher dimensional version of the Schneider--Lang Criterion,  the conclusion could be that the exceptional set of points, where all the functions take simultaneously values in a number field $K$,  is contained in an algebraic hypersurface, the bound for the number of points being replaced by a bound for the degree of the hypersurface. This program was  fulfilled by E.~Bombieri in 1970
\cite{MR0306201}. 

\begin{theorem}[Bombieri]\label{Th:Bombieri}
Let
$d$ and $n$ be two integers with $d>n\ge 1$,
$K$ be a number field, and $f_1,\ldots,f_d$ be
algebraically  independent entire functions of finite order
of growth. Assume, for $1\le \nu\le n$ and $1\le i\le d$,
that the partial derivative
$(\partial/\partial z_\nu)f_i$ of $f_i$ belongs to the ring
$K[f_1,\ldots,f_d]$. Then the set of points $w\in \bC^n$ where the $d$ functions  $f_1,\ldots,f_d$ all take values in $K$ is contained in an algebraic hypersurface. 
\end{theorem}

Bombieri produces an upper bound for the degree of such an hypersurface (see also \cite{MR88i:11047} 
 Th.~5.1.1).
His  proof \cite{MR0306201} 
 involves different tools, including \m{L^2}--estimates by {L.~Hörmander} for functions of several variables. One main difficulty that Bombieri had to overcome was to generalize Schwarz Lemma to several variables, and his solution involves an earlier work by E.~Bombieri  and S.~Lang 
\cite{MR45:5089}, 
where they use Lelong's theory of the mass of zeroes of analytic functions in a ball. Chapter 7 of 
 \cite{MR88i:11047} 
is devoted to this question. The next statement (Proposition \ref{P:Schwarz})  follows from the results in Chapter 7 of  \cite{MR88i:11047}. Given a finite subset $S$ of $\bC^n$ and an integer $t\ge 1$, one defines $\omega_t(S)$ as the smallest degree of a polynomial having a zero at each point of $S$ of multiplicity at least $t$. Then the sequence $\omega(t)/t$ has a limit $\Omega(S) $, which satisfies
$$
\frac{1}{t+n-1}\omega_t(S)\le \Omega(S)\le \frac{1}{t}\omega_t(S)\le \omega_1(S)
$$ 
for all $t\ge 1$. 

\begin{proposition}[Schwarz Lemma in several variables] \label{P:Schwarz}
Let $S$ be a finite subset of $\bC^n$ and $\epsilon$ be a positive real number. There exists a positive real number $r_0=r_0(S,\epsilon)$ such that, for any positive integer $t$, any real numbers $R$ and $r$ with $R>r\ge r_0$  and any entire function having a zero of multiplicity $\ge t$ at each point of $S$, 
$$
|f|_r\le \left(
\frac{4nr}{R}\right)^{t(\Omega(S)-\epsilon)}|f|_R.
$$
\end{proposition}

More recent results on Schwarz Lemma in several variables are due to D.~Roy \cite{MR1824984,MR1838090,MR1916273,MR1943883}.
 In particular, \cite{MR1916273} shows that Conjecture 7.1.10 of  \cite{MR88i:11047},  
on Schwarz's Lemma in several variables for finitely generated subgroups of $\bC^n$, does not hold without a technical condition of distribution. 

\subsubsection{The six exponentials Theorem}\label{SSS:SixExponentials}

Here is the six exponentials Theorem, due to Siegel, Lang and Ramachandra  \cite{MR35:5397,MR50:12931,MR1756786,MR88i:11047}.

\begin{theorem} [Six exponentials Theorem]\label{T:SixExp}
Let $x_1,\ldots,x_d$ be complex numbers which are linearly independent over $\bQ$ and let 
$y_1,\ldots,y_\ell$ be also complex numbers which are linearly independent over $\bQ$. 
Assume $d\ell>d+\ell$. Then one at least of the numbers 
\begin{equation}\label{E:SixExp}
e^{x_iy_j},\quad (1\le i\le d, \; 1\le j \le \ell)
\end{equation}
is transcendental. 
\end{theorem}

The condition $d\ell>d+\ell$ in integers $d$ and $\ell$ means that the relevant cases are $d=2$ and $\ell=3$ or $d=3$ and $\ell=2$ (there is a symmetry), hence the name of the statement  (since $\ell d=6$ in these cases). We refer to 
\cite{MR35:5397,
MR50:12931} 
for more information on this topic. 

The classical proof by Schneider's method involves an auxiliary function of the form
$$
F(z)=P(e^{x_1z},\ldots,e^{x_dz}),
$$
where the existence of the  polynomial $P$  follows from Dirichlet's box principle and the Thue--Siegel Lemma  \ref{L:TS}. The conditions which are required are that $F$ vanishes at many points of the form $s_1y_1+\ldots+s_\ell y_\ell$, with varying integers $s_1,\ldots,s_\ell$. The induction shows that $F$ vanishes at more points of this form, until one deduces that $F$ vanishes at all such points, and the conclusion then easily follows. One could avoid an infinite induction by using a zero estimate. 

A variant is to use an interpolation determinant (see \S~\ref{S:InterpolationDeterminants} and 
\cite{MR1756786} Chap.~2). 
As an example of upper bound for an interpolation determinant, here is   Lemma 2.8 in \cite{MR1756786}.

\begin{lemma}
Let $\varphi_1,\ldots,\varphi_L$ be
entire functions in $\bC$, $\zeta_1,\ldots,\zeta_L$ be elements
of $\bC$, $\sigma_1,\ldots,\sigma_L$ nonnegative integers,
and
$0<r\le R$ be real numbers, with $|\zeta_\mu|\le r$ \
($1\le\mu\le L$).  Then the absolute value of the determinant 
$$
\Delta=\det\biggl(\Bigl(
\frac{d}{ dz}\Bigr)^{\sigma_\mu}
\varphi_\lambda(\zeta_\mu)\biggr)_{1\le\lambda,\mu\le L}
$$
is bounded from above by
$$
|\Delta|\le\left(\frac{R}{
r}\right)^{-L(L-1)/2+\sigma_1+\cdots+\sigma_L}L!
\prod_{\lambda=1}^L\max_{1\le\mu\le
L}\sup_{|z|=R}\left|\left(\frac{d}{ dz}\right)^{\sigma_\mu}
\varphi_\lambda(z)\right|.
$$

\end{lemma}

Another method of proof of the six exponentials Theorem \ref{T:SixExp} is proposed in \cite{MR1462849,MR1689530}. 
The starting idea is to consider the function of two complex variables $e^{zw}$. If all numbers in (\ref{E:SixExp}) are algebraic, then this function takes algebraic values (in fact in a number field) at all points $(z,w)$, where $z$ is of the form $t_1x_1+\cdots+t_dx_d$ and $w$ of the form 
 $s_1y_1+\ldots+s_\ell y_\ell$, with integers $t_1,\ldots,t_d$, $s_1,\ldots,s_\ell$. The method does not work like this: a single function $e^{zw}$ in two variables does not suffice, one needs several functions. For this reason one introduces {\it redundant variables}.
    Letting  $z_h$ and $w_k$ be new independent  variables, one investigates the values of the functions are $e^{z_hw_k}$ at the points of Cartesian products. 
    
Redundant variables had already  been introduced in transcendental number theory in 1981 by P.~Philippon \cite{0459.10024}, for giving a proof of an algebraic independence result announced by G.V.~Chudnovskii -- at that time, the original proof 
     was very complicated  ; the approach by Philippon, using a criterion for algebraic independence due to  \'E.~Reyssat, introduced  a dramatic simplification. Now, much sharper results are known. In 1981, introducing several variables was called {\it Landau's trick}, which is a homogeneity argument: letting the number of variables tend to infinity enables one to kill error terms
      (see for instance  \cite{MR88i:11047}, \S~7.1.a p.~116). 

 The generalisation of the six exponentials Theorem to higher dimension, also with multiplicities, is one of the main topics of \cite{MR1756786}. 

Further auxiliary functions occur in the works  on algebraic independence
by  A.O.~Gel'fond \cite{MR22:2598}, 
G.V.~Chudnovskii  and others. A reference 
 is 
\cite{MR1837822}. 

\subsubsection{Modular functions}\label{SSS:Modular}

The solution, by the team at St Etienne
\cite{MR1369409},  
of the problems raised by Mahler and Manin on the transcendence of the values of the modular function $J(q)$ for algebraic values of $q$ in the unit disc,%
\footnote{The connexion between $J(q)$ and $j(\tau)$ is $j(\tau)=J(e^{2i\pi\tau})$.} %
involves an interesting auxiliary function. The general scheme of proof is the one of Gel'fond and Schneider. They construct their auxiliary function  by means of the Thue--Siegel Lemma  \ref{L:TS}  as  a polynomial in $z$ and $J(z)$  having a high multiplicity zero at the origin, say $\ge L$. They consider the exact order of multiplicity $M$, which by construction is $\ge L$, and they bound $|z^{-M} F(z)|$ in terms of $|z|$ (this amounts to using the easiest case of Schwarz Lemma in one variable with a single point). They apply this estimate to $z=q^S$, where $S$ is the smallest integer for which $F(q^S)\not=0$. Liouville's estimate gives the conclusion. See also 
\cite{MR1837822}, Chap.~2.

A variant of this construction was performed by 
Yu.~V.~Nesterenko in 1996  (see \cite{MR97m:11102} 
and 
\cite{MR1837822} Chap.~3), 
when he proved the  algebraic independence of $\pi$, $e^\pi$ and $\Gamma(1/4)$: his main result is that  {\it for any $q$ in the open set $0<|q|<1$, the transcendence degree of the field 
$$
\bQ\bigl(q,P(q),Q(q),R(q)\bigr)
$$
is at least $3$. }  Here, $P$, $Q$, $R$ are the classical Ramanujan functions, which are sometimes denoted as $E_2$, $E_4$ and $E_6$ (Eisenstein series). The auxiliary function $F$ is a polynomial in the four functions $z$, $E_2(z)$, $E_4(z)$, $E_6(z)$ ;  like for the {\it théorème stéphanois} on the transcendence of $J(q)$, it is constructed by means of the  Thue--Siegel Lemma  \ref{L:TS}, so that it has a zero at the origin of  large multiplicity, say  $M$. In order to apply a criterion for algebraic independence, Nesterenko needs to establish an upper bound for $M$, and this is not an easy result. An alternative argument, due to Philippon
(see
\cite{MR1837822} Chap.~4), 
 is to apply a measure of algebraic independence of Faisant and Philibert
\cite{MR873876} on numbers of the form $\omega/\pi$ and $\eta/\pi$.

\subsection{Universal auxiliary functions}\label{SS:UniversalAuxiliaryFunction}

\subsubsection{A general existence theorem}\label{SSS:GeneralExistence}

In the Gel'fond--Schneider method, the auxiliary function is constructed by means of the Thue--Siegel Lemma  \ref{L:TS}, and the requirement is that it has many zeroes (multiplicity are there in Gel'fond's method, not in Schneider's method). There is an alternative construction, which was initiated in a joint work with M.~Mignotte in 1974  \cite{MR0376552}, 
in connexion with quantitative statements related with transcendence criteria like the Schneider--Lang criterion. This approach turned out to be especially efficient in another context, namely in extending Schneider's method to several variables \cite{MR82k:10042}. The idea is to require that the auxiliary function $F$ has small Taylor coefficients at the origin; it follows that its modulus on some discs will be small, hence its values (including derivatives, if one wishes) at points in such a disc will also be small. Combining Liouville's estimate with Cauchy's inequality for estimating the derivatives, one deduces that $F$ has a lot of zeroes, more than would be reached by the Dirichlet's box principle. At this stage there are several situations. The easy case is when a sharp zero estimate is known: we immediately reach the conclusion without any further extrapolation: in particular there is no need of Schwarz Lemmas in several variables. This is what happend in 
\cite{MR82k:10042} 
for exponential functions in several variables, the zero estimate being due to D.W.~Masser
\cite{MR82j:10062}. 
This transcendence result, dealing with products of multiplicative groups (tori), can be extended to commutative algebraic groups 
\cite{MR84j:10046}, 
thanks to the zero estimate of Masser and Wüstholz
\cite{MR83d:10040,MR86h:11054}.

A {\it zero estimate} is a statement which gives a lower bound for the degree of a polynomial which vanishes at given points (when multiplicities are introduced, this is sometimes called a {\it multiplicity estimate}). D.W.~Masser developed the study of zero estimates, and P.~Philippon, G.~Wüstholz, Yu.V.~Nesterenko were among those who contributed to the theory.
Introductions to this subject have been written by D.~Roy (Chapters 5 and 8 of 
\cite{MR1756786} and Chapter 11 of \cite{MR1837822}). Also Chapter 5 of \cite{pre05232987} is devoted to multiplicity estimates.

A {\it Schwarz Lemma} is a statement which gives an upper bound for the maximum modulus of an analytic function which vanishes at given points -- multiplicities may be there. We gave an example in Proposition \ref{P:Schwarz}.

When the function has only small values at these points, instead of zeroes, one speaks either of a {\it small value Lemma} or of an {\it approximate Schwarz Lemma}
\cite{MR1453214}. 

A Schwarz Lemma implies a zero estimate, and a small value Lemma implies a Schwarz Lemma. However, the assumptions for obtaining  a small value Lemma, for instance, are usually stronger than for only a Schwarz Lemma: as an example, in one variable, there is no need to introduce an assumption on the distance of the given points for a Schwarz Lemma, while  a small value Lemma woud not be true without such a condition.

The construction of universal auxiliary functions  is developed in  
\cite{MR1243105} (see Lemme 2.1) 
and
\cite{MR1176537}. 
Here is Proposition 4.10 of 
\cite{MR1756786}.

\begin{proposition}\label{P:UniversalAuxiliaryFunction}
Let $L$ and $n$ be positive
integers, $N$, $U$, $V$, $R$, $r$ positive real numbers and
$\varphi_1,\ldots,\varphi_L$ entire functions in $\bC^n$.
Define $W=N+U+V$ and assume 
$$
W\ge 12 n^2,\quad  e\le \frac{R}{ r}\le e^{W/6},\quad
\sum_{\lambda=1}^L|\varphi_\lambda|_R\le e^U
$$
and
$$
(2W)^{n+1}\le LN\bigl(\log(R/r)\bigr)^n.
$$
Then there exist rational integers $p_1,\ldots,p_L$, with
$$
0<\max_{1\le\lambda\le L}|p_\lambda|\le e^N,
$$
such that the function $F=p_1\varphi_1+\cdots+p_L\varphi_L$
satisfies 
$$
|F|_r\le e^{-V}.
$$ 
\end{proposition}

For an application to algebraic independence, see 
 \cite{MR1837822} Prop. 3.3 Ch.~14. 

 The construction of universal auxiliary functions is one of the tools in D.~Roy's approach to Schanuel's Conjecture in \cite{MR1824984,MR1838090}.

\subsubsection{Duality}\label{SSS:Duality}

In the papers \cite{MR1243105} (see Lemme 2.1) 
and
\cite{MR1176537},  a {\it dual} construction is performed, where auxiliary {\it analytic functionals} are constructed. The duality 
(\ref{Eq:DualiteGS}) between the methods of Schneider and Gel'fond  can be  explained by means of the Fourier--Borel transform
(see \cite{MR1243105}, especially Lemmas 3.1 and 7.6, 
and  \cite{MR1756786},  \S~13.7). 
 In the special case of exponential polynomials,  the Fourier--Borel duality reduces to the relation 
{
\def\utau{{\underline{\tau}}}
\def\ut{{\underline{t}}}
\def\usigma{{\underline{\sigma}}}
\def\uw{\underline{w}}
\begin{equation}\label{E:duality}
D^{\usigma} 
\left(
\uz^{\utau} e^{\ut\uz}
\right)(\us)=
D^{\utau} 
\left(
\uz^{\usigma}
e^{\us\uz}
\right)
(\ut),
\end{equation}
where $\usigma$,  $\utau$, $\us$, $\ut$,  $\uz$   stand  for tuples
$$
\begin{array}{ll}
\usigma=(\sigma_1,\ldots,\sigma_n)\in\bZ_{\ge 0}^n,
&
\utau=(\tau_1,\ldots,\tau_n)\in\bZ_{\ge 0}^n,
\\
\us=(s_1,\ldots,s_n)\in\bC^n,
&
\ut=(t_1,\ldots,t_n)\in\bC^n,
\\
\uz=(z_1,\ldots,z_n)\in\bC^n
&
\end{array}
 $$
and
$$
 D^{\usigma} =\left(\frac{\partial}{\partial z_1}\right)^{\sigma_1}
 \cdots
 \left(\frac{\partial}{\partial z_n}\right)^{\sigma_n}
$$  
 (see   \cite{MR1243105} Lemma 3.1 and \cite{MR1756786} Corollary 13.21). This is a generalization in several variables of the formula
 $$
 \left( \frac{d}{dz}\right)^{\sigma}\left(
z^{\tau} e^{tz}\right)
 (s)=
 \left( \frac{d}{dz}\right)^{\tau}\left(
z^{\sigma}
e^{sz}
\right)
(t)
$$
for $t$, $s$  in $\bC$ and $\sigma$, $\tau$ non--negative integers, both sides being
$$
\sum_{k=0}^{\min
\{
\tau,\sigma
\}
}
\frac{\sigma!\tau!}
{k! (\tau-k)!(\sigma-k)!}
t^{\sigma-k} s^{\tau-k}e^{st}.
$$  

Another (less special) case of (\ref{E:duality}) is the duality between  Schneider's method in several variables and Baker's method.  In Baker's method, one considers the values at
$ ( t, t\log\alpha_1,\ldots,t\log\alpha_n)$ of 
$$
 \left(\frac{\partial}{\partial z_0}\right)^{\tau_0}
  \cdots 
   \left(\frac{\partial}{\partial z_n}\right)^{\tau_n}
    \bigl( 
    z_0^{\sigma}e^{s_1z_1}\cdots e^{s_nz_n}
    e^{s_0(\beta_0z_0+\cdots+\beta_nz_n)}
    \bigr),
$$
while Schneider's method in several variables deals with the values of 
$$
 \left(\frac{\partial}{\partial z_0}\right)^{\sigma}
 \bigl( 
 z_0^{\tau_0} z_1^{\tau_1} \cdots  z_n^{\tau_n}
 (e^{z_0}\alpha_1^{z_1}\cdots \alpha_n^{z_n})^t 
 \bigr)
 $$
 at the points 
$$
 (s_0\beta_0,s_1+s_0\beta_1, \ldots, s_n+s_0\beta_n).
$$
Again, these values are just the same, namely
$$
\sum_{k=0}^{\min
\{
\tau_0,\sigma
\}
}
\frac{\sigma!\tau_0!}
{k! (\tau_0-k)!(\sigma-k)!}
t^{\sigma-k} (s_0\beta_0)^{\tau_0-k}
(s_1+s_0\beta_1)^{\tau_1}\cdots
(s_n+s_0\beta_n)^{\tau_n} 
(\alpha_0^{s_0}\cdots \alpha_n^{s_n})^t
$$ 
with 
$$
\alpha_0=e^{\beta_0}\alpha_1^{\beta_1}\cdots\alpha_n^{\beta_n}.
$$
This is a special case of (\ref{E:duality}) with 
$$
\utau=(\tau_0,\ldots,\tau_n)\in\bZ_{\ge 0}^{n+1},
\quad
\usigma=(\sigma_0,0,\ldots,0)\in\bZ_{\ge 0}^{n+1},
$$
$$
\ut=(t,t\log\alpha_1,\ldots,t\log\alpha_n)\in\bC^{n+1},
\quad 
\us=(s_0\beta_0, s_1+s_0\beta_1,\ldots,s_n+s_0\beta_n)\in\bC^{n+1},
$$
$$ 
\uz=(z_0,\ldots,z_n)\in\bC^{n+1}.
$$ 
}

\subsection{Mahler's Method}\label{SS:MahlerMethod}

In 1929, K.~Mahler
\cite{JFM55.0115.01}
 developed an original method to prove the transcendence of values of functions satisfying certain types of functional equations. This method was somehow forgotten, for instance it is not quoted among  the 440 references of the survey paper
 \cite{MR0214551} 
 by N.I.~Fel'dman, and A.B.~{\v{S}}idlovski{\u\i}.  After the publication of the paper
 \cite{MR40:2611} 
 by Mahler in 1969, several mathematicians 
(including K.~Kubota, J.H.~Loxton and  A.J.~van der Poorten) extended the method
 (see the Lecture Notes \cite{MR1439966} 
 by K.~Nishioka for further references).  
 The construction of the auxiliary function is similar to what is done in Gel'fond--Schneider's method, with a main difference: in place of the Thue--Siegel Lemma  \ref{L:TS}, only linear algebra is required. No estimate for the coefficients of the auxiliary polynomial is needed   in Mahler's method.

 The following example is taken from \S~1.1 of  \cite{MR1439966}. Let $d\ge 2$ be a rational integer and $\alpha$ an algebraic number with $0<|\alpha|<1$. Here is a sketch of proof that the number
 $$
 \sum_{k=0}^\infty \alpha^{d^k}
 $$
 is transcendental, a result due to K.~Mahler \cite{JFM55.0115.01}. The basic remark is that the function
 $$
f(z)= \sum_{k=0}^\infty z^{d^k}
 $$
 satisfies the functional equation $f(z^d)=f(z)-z$. It is not difficult to check also that $f$ is a transcendental function, which means that if $P$ is a non--zero polynomial in $\bC[X,Y]$, then the function $F(z)=P\bigl(z,f(z)\bigr)$ is not the zero function. From linear algebra, it follows that if $L$ is a sufficiently large integer,  there exist a non--zero polynomial $P$  in $\bZ[X,Y]$ of partial degrees $\le L$ such that the associated function  $F(z)=P\bigl(z,f(z)\bigr)$ has a zero at the origin of multiplicity $> L^2$. Indeed, the existence of $P$ amounts to showing the existence of a non--trivial solution to a system of $L^2+1$ homogeneous linear equations with rational coefficients in $(L+1)^2$ unknowns.

If $T$ denotes the multiplicity of $F$ at the origin, then the limit
 $$
 \lim_{k\rightarrow\infty} F(\alpha^{d^k})\alpha^{-d^kT}
 $$
 is a non--zero constant. One deduces an upper bound for $|F(\alpha^{d^k})|$ when $k$ is a sufficiently large  positive integer, and this upper bound is not compatible with Liouville's inequality. Hence $f(\alpha)$ is transcendental. 
 
 Another instructive example of application of Mahler's method is given by D.W.~Masser in the first of his Cetraro's lectures
 \cite{MR2009828}. 
 
  The numbers whose transcendence is proved by this method are not Liouville numbers, but they are quite well approximated by algebraic numbers. In the example we discussed, the number $f(\alpha)$ is very well approximated by the algebraic numbers
   $$
 \sum_{k=0}^K \alpha^{d^k}, \qquad (K>0).
 $$

\section{Interpolation determinants}\label{S:InterpolationDeterminants}

An interesting development of the saga of auxiliary functions took place in 1991, with the introduction of interpolation determinants by M.~Laurent. Its origin goes back to earlier works on a question raised by Lehmer, which we first discuss. 

\subsection{Lehmer's Problem}

Let \m{\theta} be a non--zero algebraic integer of degree \m{d}. Mahler's {\it measure } of \m{\theta} is  
\M{
M(\theta)=\prod_{i=1}^d\max(1,|\theta_i|)
=\exp\left(
\int_0^1 \log|f(e^{2i\pi t}| dt\right),
}
 where \m{\theta=\theta_1} and \m{\theta_2,\cdots,\theta_d} are the conjugates of \m{\theta} and \m{f} the (monic) irreducible polynomial of \m{\theta} in \m{\bZ[X]}
 $$
 f(X)=(X-\theta_1)\cdots(X-\theta_d)\in\bZ[X].
 $$
From the definition one deduces \m{M(\theta)\ge 1}. According to a well-known and easy result of Kronecker,    \m{M(\theta)=1} {\it  if and only if \m{\theta}   is a root of unity. }

D.H.~Lehmer
\cite{Lehmer}
 asked whether {\it  there is a constant \m{c>1} such that \m{M(\theta)<c} implies that \m{\theta} is a root of unity. }

Among many tools which have been introduced to answer this question, we only quote some of them which are relevant for our concern. 
In 1977, M. Mignotte
\cite{MR80b:12002}
 used ordinary Vandermonde determinants to study algebraic numbers whose conjugates are close to the unit circle. 
 In 1978, C.L.~Stewart
\cite{MR507748} 
sharpened earlier results by Schinzel and Zassenhaus (1965) and 
Blanksby and Montgomery (1971) by
introducing an auxiliary function (whose existence follows from the Thue--Siegel Lemma \ref{L:TS}) and using an extrapolation similar to what is done in the Gel'fond--Schneider method.

Auxiliary functions, introduced by Stewart, were one of the tools used by  
 E.~Dobrowolski in 1979  
 \cite{MR543210} 
His ability to exploit congruences mod $p$ was a major advance that significantly sharpened the previous estimates of Blanksby--Montgomery  
and Stewart
.  This use of congruences
mod $p$, and then summing the estimates over $p$ in a suitable interval, resulted in a very substantial
improvement, falling just short of the Lehmer conjecture.  He achieved the best unconditional result known so far in this direction  (apart from some marginal improvements on the numerical value for \m{c}): 
\begin{quote}
 {\it There is a constant \m{c} such that, for \m{\theta} a non--zero algebraic integer of degree \m{d},
 \M{
 M(\theta)<1+c(\log\text{}\log d/\log d)^3
 }
 implies that \m{\theta} is a root of unity. 
 }
\end{quote}
In 1982,  D.~Cantor and E.G.~Straus
\cite{MR679001} 
revisited this method of Dobrowolski by replacing the  auxiliary function by a generalised Vandermonde determinant. The idea is the following: in Dobrowolski's proof, there is a zero lemma which can be translated into a statement that some matrix has a maximal rank; therefore, some determinant is not zero. On the one hand, this determinant is not too large: an upper bound for its absolute value  follows from Hadamard's inequality; the upper bound depends on $M(\theta)$. On the other hand, the absolute value of this determinant is shown to be big, because it has many factors of the  form \m{\prod_{i,j}|\theta_i^p-\theta_j|^k}, for many primes \m{p}. The lower bound makes use of a Lemma due to Dobrowolski:  {\it For  \m{\theta} not a root of unity, 
\M{
\prod_{i,j}|\theta_i^p-\theta_j|\geq p^d
} 
for any prime \m{p}.
}
Combining the upper and lower bounds yields the conclusion. 

One may also prove the lower bound by means of a $p$--adic Schwarz Lemma: a function (here merely a polynomial) with many zeroes has a small ($p$--adic) absolute value. 
This alternative argument produces a proof which is similar to the  ones arising from transcendental number theory, with analytic estimates on the one side and arithmetic estimates (Liouville type, or product formula) on the other. See \cite{MR1756786} \S~3.6.5 and  \S~3.6.6.
 
Dobrowolski's result has been extended to several variables by F.~Amoroso and S.~David in 
\cite{1011.11045} -- the higher dimensional version is much more involved.  We refer to S.~David's survey \cite{SinnouRamanujanJM}  for further references on this topic. 
We only notice the role of the generalized Dirichlet exponents $\omega_t$ and $\Omega$ (see  \cite{MR88i:11047}, \S~1.3 and Chap.~7), 
which were introduced in transcendental number theory, in connexion with multidimensional Schwarz Lemmas (cf{.} Proposition \ref{P:Schwarz}).

\subsection{Laurent's interpolation determinants}

In 1991, M.~Laurent 
\cite{MR1144324} 
discovered that one may get rid of the Dirichlet's box principle in Gel'fond--Schneider's method by means of his {\it interpolation determinants}. In the classical approach, there is a zero estimate (or vanishing estimate, also called multiplicity estimate when derivatives are there), which shows that some auxiliary function cannot have too many zeroes. This statement can be converted into the non--vanishing of some determinant. Laurent works directly with this determinant: a Liouville-type estimate produces a lower bound; the remarkable fact is that analytic estimates like Schwarz Lemma produce sharp upper bounds. Again, the analytic estimates depend only in one variable (even if the determinant is a value of a function in many variables, it suffices to restrict this function to a complex line -- 
see \cite{MR1756786} 
\S~6.2). Therefore, this approach is especially efficient when dealing with functions of several variables, where Schwarz Lemmas are lacking. 

Baker's method is introduced in \cite{MR1756786} \S~10.2 with interpolation determinants 
and in \S~10.3 with auxiliary functions.

Interpolation determinants are easy to use when a sharp zero estimate is available. If not, it is more tricky to prove the analytic estimate. However, it is possible to perform extrapolation in the transcendence methods involving interpolation determinants: an example is a proof of 
P\'olya's theorem 
(\ref{E:Polya}) 
by means of interpolation determinants 
\cite{MR1659778}.  

Proving algebraic independence results by means of interpolation determinants was quoted as an open problem in \cite{MR1243105} p.~257: at that time (1991), auxiliary functions (or auxiliary functionals) were required, together with a zero estimate (or an interpolation estimate). The first  solution to this problem in 1997
\cite{MR1476295} 
(see also \cite{MR1756786} Corollary 15.10)
makes a detour via measures of simultaneous approximations: such measures can be proved by means of interpolation determinants, and they suffice to produce algebraic independence statements   ({\it small transcendence degree}: at least two numbers in certain given sets are algebraically independent). 
Conjecturally, this method should  produce {\it large transcendence degree} results as well  -- see \cite{MR1756786} Conjecture 15.31. See also the work by P.~Philippon 
\cite{MR1752252}.

Another approach, due to M.~Laurent and D.~Roy \cite{MR1837427} (see also \cite{MR1837822} Chap.~13), 
is based on the observation that in algebraic independence proofs, the determinants which occur produce sequences of polynomials having small values {\it together with their derivatives} at a given point. By means of a generalization of Gel'fond's transcendence criterion involving  multiplicities, 
M.~Laurent and D.~Roy succeeded to get algebraic independence results. Further generalisations of Gel'fond's criterion involving not only multiplicities, but also several points (having some structure, either additive or multiplicative) are being investigated by D.~Roy in 
connection with his original strategy towards a proof  of  Schanuel's Conjecture
\cite{MR1824984,
MR1838090
}.

\section{Bost slope inequalities, Arakelov's Theory}\label{S:Arakelov}

Interpolation determinants require choices of bases. A further tool has been introduced by
J-B.~Bost in 1994 \cite{MR1423622}, 
where bases are no longer required: the method is more intrinsic. His argument  rests on Arakelov's Theory, which is used to  produce
 {\it slope inequalities}.  This new approach is especially interesting for results on abelian varieties obtained by transcendence methods, the examples developed by Bost being related with the work of D. Masser and G. W\"ustholz on periods and isogenies of abelian varieties over number fields. Further estimates related to Baker's method and measures of linear independence of logarithms of algebraic points on abelian varieties have been achieved by E.~Gaudron \cite{MR2003a:11090,MR2292632,MR2362434} 
using Bost's approach. An enlightening introduction to Bost method can be found in the Bourbaki lecture \cite{MR1975179}  by 
A.~Chambert-Loir  
in 2002.

\def\cprime{$'$} \def\cprime{$'$} \def\cprime{$'$} \def\cprime{$'$}
  \def\cprime{$'$}

\vfill

 \vskip 2truecm plus .5truecm minus .5truecm 

\font\ninerm=cmr10 at 9pt

\font\ninett=cmtt10 scaled 900 

\hfill
\vbox{\ninerm
 \hbox{Michel WALDSCHMIDT}
 \hbox{Universit\'e Pierre et Marie Curie -- Paris 6}
 \hbox{IMJ UMR 7586} 
 \hbox{175, rue du Chevaleret }
 \hbox{PARIS F--75013 FRANCE}
 \hbox{e-mail: 
\href{mailto:miw@math.jussieu.fr}{miw@math.jussieu.fr}}
 \hbox{URL:
\href{http://www.math.jussieu.fr/~miw/}%
{http{$:$}//www.math.jussieu.fr/${\scriptstyle \sim}$miw/}}
} 

\vfill

  \end{document}